\documentclass[12pt]{amsart}
\usepackage{amsmath,amssymb,latexsym,amsthm, amscd}
\usepackage{mathrsfs}
\theoremstyle{definition}
\theoremstyle{plain}

\newtheorem{thm}[subsection]{Theorem}
\newtheorem{lem}[subsection]{Lemma}
\newtheorem{cor}[subsection]{Corollary}

\newcommand{\mbf}{\mathbf}
\newcommand{\mrm}{\mathrm}
\newcommand{\bin}[2]{
\left[
   \begin{array}{@{}c@{}}
     #1 \\#2
   \end{array}
\right]    }
 
\title[A realization of quantum groups ]
{A realization of quantum groups via product valued quivers} 
\author{ Yiqiang Li}
\address{Department of Mathematics\\ Kansas State University\\ Manhattan, KS  66506}
\email{yqli@math.ksu.edu}
\author{Zongzhu Lin}
\address{Department of Mathematics\\ Kansas State University\\ Manhattan, KS  66506}
\email{zlin@math.ksu.edu}
\thanks{The research was partially supported by NSF grant DMS-0200673.\\ The main results in this paper 
were reported by Y. Li  in the Joint Mathematics Meetings at San Anotonio, Texas, January 12-15, 2006.  
They  were also reported by Z. Lin in the algebra seminar at Kansas State University, April 10, 2006.}
\begin{document}

\begin{abstract} 
Let $\overset{\rightarrow}{\Gamma}$ be a valued quiver.
Let $\mrm C$ be the symmetrizable generalized Cartan matrix associated
to $\overset{\rightarrow}{\Gamma}$. We show that 
the whole quantum group associated to $\mrm C$ can be realized from 
the category of the representations of the product valued quiver 
$\overset{\rightarrow}{\Gamma^{\pm}}$.
This method can be used to realize the whole generalized Kac-Moody
Lie algebra associated to $\mrm C$, as discussed in ~\cite{Li-Lin}.
\end{abstract}

\maketitle

\section{Introduction}
Let $\Lambda$ be a finite dimensional hereditary algebra 
over a finite field $\mbf F_q$.
Let $\mrm C$ be the associated symmetrizable Cartan matrix.
It is well-known that the composition algebra of $\Lambda$ gives a  
realization of the positive part $\mbf U(\mrm C)^-$ of 
the quantized enveloping algebra associated to $\mrm C$  
(specialized at $v=q^{\frac{1}{2}}$), thanks to ~\cite{Ringel2} and ~\cite{Green}.
As Ringel mentioned in ~\cite{Ringel2}, 
it is of interest to recover the whole quantized enveloping algebra,
instead of just positive part,
from the module categories of finite dimensional hereditary algebras.
Xiao gave an answer in \cite{Xiao} by using the Drinfeld double 
to piece together two Borel parts.
Later on, several attempts have been made to recover the whole quantum
group from different settings. (See ~\cite{Kapranov2} and ~\cite{Toen}.)

In this paper, we give a new realization of the whole quantum group
from the module category of certain finite dimensional hereditary
algebra. 
Instead of considering the category, $\Lambda$-mod,  
of finite dimensional $\Lambda$-modules, 
we consider the category of $\Lambda^{\pm}$-modules,
which can be thought of as the double of $\Lambda$-mod.
(See Section 3.)
By tensoring the torus part with the twisted Ringel-Hall algebra 
$\mbf H(\Lambda^{\pm})$ and taking a quotient modulo  some relations, 
we recover the whole quantum group $\mbf U(\mrm C)$ (specialized at $v=q^{\frac{1}{2}}$).
The double Ringel-Hall algebra (\cite{DX1}) associated to $\Lambda$
is also recovered along the way of the construction.
Moreover, this construction can be generalized to the case 
when the valued quiver $\overset{\rightarrow}{\Gamma^{\pm}}$ is replaced by 
the  valued quiver $\overrightarrow{\Gamma}^{2n}$ 
($n$ a positive integer) in Section 4. 

In the forthcoming paper ~\cite{Li-Lin}, the realization of the whole Lie algebras from
the categories of representations of locally finite quivers will
be further discussed.
The relation with the affinization of semisimple Lie algebras will be
addressed.

The paper is organized as follows.
In Section 2, we give the construction on the theoretical level. 
In Section 3,
we apply the construction in Section 2 to the module categories of 
finite dimensional hereditary algebras, associated to a product valued
quiver. 
The quantized enveloping algebra $\mbf U(\mrm C)$ associated to $\mrm C$ is realized.
Section 4 generalizes the previous constructions. 

{\bf Acknowledgements.} We wish to thank Professor Bangming Deng for
helpful discussions.

\section{Quantum groups}

\subsection{Definitions.}
Let $I$ be an index set.
Let $\mrm{C}=(a_{ij})_{i,j\in I}$ be 
a symmetrizable generalized Cartan matrix. 
There are integers $d_i$ with no proper common divisor such that 
$d_i a_{ij}=d_j a_{ji}$ for all $i, j $ in $I$.
Let $v$ be an indeterminate. As in [L3], we set
\begin{itemize}
\item [(1)] $v_i=v^{d_i}\quad$  for all $i\in I$.
\item [(2)] $[m]=\frac{v^m-v^{-m}}{v-v^{-1}}$ and
           $[m]^!=\prod^{m}_{t=1}[t]$ $\quad$ for any $m\in \mathbb N$.
\item [(3)] $\bin{m}{t}= \frac{[m]^!}{[t]^![m-t]^!}$,
           and 
           $[m]_i=\frac{v_i^m-v_i^{-m}}{v_i-v_i^{-1}}$ $\quad$ for any $m, t\in \mathbb N$ such
           that $0\leq t\leq m$.
\end{itemize}

The quantized enveloping algebra $\mbf U=\mbf U(\mrm C)$ associated to $\mrm C$ is 
a $\mathbb Q(v)$-algebra characterized by generators and relations.
The generators are
\[
E_i, F_i, K_i \; \text{ and} \; K_{-i} \quad   \text{for all}\; i \in I
\]
and the relations are
\begin{itemize}
   \item[($\mbf U1$)] $K_0=1$, $K_iK_{-i}=1=K_{-i}K_i$, and $K_iK_j=K_jK_i\quad$ 
             for any $i,j \in I$;
   \item[($\mbf U2$)] $K_iE_j=v_i^{a_{ij}}E_jK_i$, $K_iF_j=v_i^{-a_{ij}}F_jK_i\quad$
             for any $i,j\in I$;
   \item[($\mbf U3$)] $E_iF_i-F_iE_i=\frac{K_i-K_{-i}}{v_i-v_i^{-1}}$  and $E_iF_j=F_jE_i\quad$
             for any $i\neq j \in I$;
   \item[($\mbf U4$)] $\sum_{t=0}^{1-a_{ij}}(-1)^t \bin{1-a_{ij}}{t}_i
             E_i^t E_j E_i^{1-a_{ij}-t}=0\quad$
             for  $i\neq j \in I$;
   \item[($\mbf U5$)] $\sum_{t=0}^{1-a_{ij}}(-1)^t \bin{1-a_{ij}}{t}_i
             F_i^t F_j F_i^{1-a_{ij}-t}=0\quad$
             for $i\neq j \in I$.
\end{itemize}

Given any $\mrm C$,
define the matrix 
$\mrm{C}^{\pm}=(A_{(\alpha,i),(\beta,j})$, 
where $(\alpha,i),(\beta,j)\in \{+,-\}\times I$ as follows.
(For convenience, we write $A_{(\alpha,i),(\beta,j)}=A^{\alpha,\beta}_{i,j}$.)
\begin{itemize}
\item $A^{+,+}_{i,j}=a_{ij}$ for $i,j\in I$;
\item $A^{-,-}_{i,j}=a_{ij}$ for $i,j\in I$;
\item $A^{+,-}_{i,i}=-2$ and $A^{-,+}_{i,i}=-2$ for $i\in I$;
\item $A^{+,-}_{i,j}=0$ and $A^{-,+}_{i,j}=0$ for $i\neq j$.
\end{itemize}
Pictorially, $\mrm{C}^{\pm}$ is the matrix
\[
 \left[
  \begin{array}{cc}
    \mrm{C} & -2\, \mrm{Id}\\
  -2\, \mrm{Id} &   \mrm{C}
  \end{array}
 \right],
\]
where $\mrm{Id}$ stands for the identity $|I|\times |I|$ matrix. 
One can check that $\mrm{C}^{\pm}$ is a symmetrizable generalized Cartan matrix.
Let $d^+_{i}=d_i$ and $d^-_i=d_i$ for all $i\in I$. Then the 
set $\{d^{\alpha}_i \;|\; (\alpha,i)\in \{+,-\}\times I\}$ is the minimal symmetrization 
of the Cartan matrix $\mrm{C}^{\pm}$.

\subsection{The $\mathfrak{sl}_2$ case.}
For example, when $\mrm{C}=[2]$, the $1\times 1$ matrix with entry 2, $\mrm{C}^{\pm}$ is the 
matrix
\[ 
 \left[
  \begin{array}{cc}
    2  & -2\\
   -2  &  2
  \end{array}
 \right].
\]
Let $\mbf{U}^+(C^{\pm})$ 
be the positive part of the quantum group associated to $C^{\pm}$ in this example. 
It is a $\mathbb{Q}(v)$-algebra characterized by generators and relations, 
with the generators:
\[
 E^+ \;\text{and}\; E^{-}
\]
and the relations:
\begin{itemize}
  \item [(a)] $(E^{+})^3E^--(v^2+1+v^{-2})(E^+)^2E^-E^+
                    +(v^2+1+v^{-2})E^{+}E^{-}(E^+)^2
                    -E^{-}(E^+)^3=0;$
  \item [(b)] $(E^-)^3E^+-(v^2+1+v^{-2})(E^-)^2E^+E^-+
                    (v^2+1+v^{-2})E^-E^+(E^-)^2
                    -E^{+}(E^-)^3=0.$
\end{itemize}

Consider the $\mathbb{Q}(v)$-vector space 
$\mathcal{V}'=\mathbb{Q}(v)[K, K^{-1}]\otimes  \mbf U^+(\mrm C^{\pm})$. 
We write $fx=f\otimes x$, 
for any $f\in \mathbb{Q}(v)[K,K^{-1}]$ and $x\in \mbf{U}^+(\mrm C^{\pm})$.
There is a unique associative $\mathbb{Q}(v)$-algebra structure on $\mathcal{V}'$ such that 
\begin{enumerate}
\item[(1)] $KK^{-1}=1=K^{-1}K$;
\item[(2)] $KE^+=v^2E^+K$;
\item[(3)] $KE^-=v^{-2}E^-K$
\end{enumerate} 
and preserves the $\mathbb{Q}(v)$-algebra structure on $\mbf U^+(\mrm C^{\pm})$. 

Define the algebra $\mathcal V$ to be the quotient of the algebra 
$\mathcal V'$ by the two-sided ideal generated by the following relation.
\begin{itemize}
\item[(4).] $E^+E^--E^-E^+=\frac{K-K^{-1}}{v-v^{-1}}$.
\end{itemize}
Then, $\mathcal V$ is a $\mathbb{Q}(v)$-algebra 
defined by the  following generators and relations.
\begin{itemize}
\item Generators: $E^+$, $E^-$, $K$ and $K^{-1}$.
\item Relations: (1), (2), (3), (4), (a), and (b) in this section.
\end{itemize}
We have 
\begin{lem}
\label{sl2}
There is a unique $\mathbb{Q}(v)$-algebra isomorphism 
$\varphi: \mbf U(\mathfrak{sl}_2) \to \mathcal V$ such that 
\[
\varphi(E)=E^+, \quad
\varphi(F)=E^-, \quad
\varphi(K)=K \quad \text{and} \quad \varphi(K^{-1})=K^{-1},
\] 
where $\mbf U(\mathfrak{sl}_2)$ is the quantized enveloping algebra 
associated to the Cartan matrix $\mrm{C}=[2]$.
\end{lem}
\begin{proof}
Note that the quantum group $\mbf U(\mathfrak{sl}_2)$ 
associated to the Cartan matrix $\mrm{C}=[2]$
is defined by the following generators and relations. 
\begin{itemize}
\item Generators: $E$, $F$, $K$ and $K^{-1}$.
\item Relations: (1), (2), (3) and (4) with $\{E^+,E^-\}$ replaced by $\{E,F\}$.
\end{itemize}
By the universality of $\mbf U(\mathfrak{sl}_2)$, 
there is a unique $\mbf Q(v)$-algebra homomorphism
$\varphi: \mbf U(\mathfrak{sl}_2) \to \mathcal{V}$ such that $\varphi(E)=E^+$,
$\varphi(F)=E^-$, $\varphi(K)=K$ and $\varphi(K^{-1})=K^{-1}$.
$\varphi$ is surjective since $E^+, E^-, K, K^{-1}$ are generators of $\mathcal V$. 

To prove that $\varphi$ is an isomorphism,
it suffices to show that in the defining relations of $\mathcal{V}$, 
the relations (1)-(4) imply the relations (a) and (b).
(Then the assignments $E^+ \mapsto E$, $E^-\mapsto F$, $K\mapsto K$ and $K^{-1}\mapsto K^{-1}$ 
define a $\mathbb Q(v)$-algebra homomorphism $\psi: \mathcal V \to \mathbf U(\mathfrak{sl}_2)$ 
and $\psi$ is the  inverse of the morphism  $\varphi$.)
Now write the left hand side of (a) in the following form.
\begin{itemize}
\item [($a'$)]
          $(E^+)^2(E^+E^--E^-E^+)-(v^2+v^{-2})E^+(E^+E^--E^-E^+)E^++
           (E^+E^--E^-E^+)(E^+)^2.$
\end{itemize}
By (4), ($a'$) is equal to 
\begin{itemize}
\item [($a''$)]
   $\frac{1}{v-v^{-1}}\left((E^+)^2(K-K^{-1})-(v^2+v^{-2})E^+(K-K^{-1})E^++
   (K-K^{-1})(E^+)^2)\right).$
\end{itemize}
By (2), we have 
\begin{itemize}
\item [($2'$).]
$(E^+)^2(K-K^{-1})=v^{-2}E^+KE^+-v^2E^+K^{-1}E^{+},$
\end{itemize}
and 
\begin{itemize}
\item[($2''$).]
$(K-K^{-1})(E^+)^2=v^2E^+KE^+-v^{-2}E^+K^{-1}E^+.$
\end{itemize}
It is obvious that ($2'$) and ($2''$) imply that ($a''$) is zero. 
Therefore, (2) and (4) imply (a). Similarly, (3) and (4) imply (b). 
\end{proof}

\subsection{General cases.}
Now given any symmetrizable generalized Cartan matrix $\mrm{C}$. Consider the positive
part $\mbf{U}^+(\mrm{C}^{\pm})$ 
of the quantum group associated to $\mrm{C}^{\pm}$. 
It is a $\mathbb{Q}(v)$-algebra defined by the following generators and relations.
The generators are
\[
E^{\alpha}_i\quad\quad \text{ where} \; (\alpha,i)\in \{+,-\}\times I 
\] 
and the relations are 
 \[
 \sum_{t=0}^{1-A^{\alpha, \beta}_{i,j}}(-1)^t
 \bin{ 1-A^{\alpha,\beta}_{i,j}}{t}_{(\alpha,i)}
 (E^{\alpha}_i)^t E^{\beta}_j (E^{\alpha}_i)^{1-A^{\alpha,\beta}_{i,j}-t}=0,
\]
where $(\alpha, i)\neq (\beta,j)\in \{+,-\}\times I$.

More explicitly, it is a $\mathbb{Q}(v)$-algebra defined by the following generators 
and relations.
The generators are
\[
E^+_{i}\;\text{ and}\; E^{-}_{i}\quad \text{ for all}\; i\in I
\]
and the relations are
\begin{align}
&\sum^{1-a_{ij}}_{t=0}(-1)^t \bin{1-a_{ij}}{t}_i
(E^+_i)^t E^+_j (E^+_i)^{1-a_{ij}-t}=0\quad 
\text{ for} i\neq j \in I.
\tag{+ +}\\  
&\tag{ - \ -}
\sum^{1-a_{ij}}_{t=0}(-1)^t \bin{1-a_{ij}}{t}_i
(E^-_i)^t E^-_j (E^-_i)^{1-a_{ij}-t}=0\quad
\text{for}\; i\neq j \in I.\\
&\tag {(+\ -)(a)} 
(E^+_i)^3E^-_i-a(E^+_i)^2E^-_iE^+_i+ aE^+_iE^-_i(E^+_i)^2-E^-_i(E^+_i)^3=0.\\
&\tag{(+\ -)(b)}
(E^-_i)^3E^+_i-a(E^-_i)^2E^+_iE^-_i+ aE^-_iE^+_i(E^-_i)^2-E^+_i(E^-_i)^3=0.\\
&\tag{(+\ -)(c)}
E^+_iE^-_j=E^-_jE^+_i\quad\text{for}\; i\neq j\in I.
\end{align}
Here $a=v^2_i+1+v^{-2}_i$.

Now consider the $\mathbb{Q}(v)$-vector space
$\mathcal{V}'=\mathbb{Q}(v)[K_{i},K_{-i}]_{i\in I}
\otimes\mbf{U}^+(\mrm C^{\pm})$.
There is a unique associative $\mathbb{Q}(v)$-algebra structure on $\mathcal{V}'$ such that 
\begin{enumerate}
\item[($\mathcal{V}1$)] $K_0=1$, $K_iK_{-i}=1=K_{-i}K_i$ and $K_iK_j=K_{i+j}\quad$ 
                       for any $i,j \in I$,
\item[($\mathcal{V}2$)] $K_iE^+_j=v_i^{a_{ij}}E^+_jK_i$ and  $K_iE^-_j=v_i^{-a_{ij}}E^-_jK_i\quad $
                       for any $i,j\in I$,
\end{enumerate} 
and preserves the $\mathbb{Q}(v)$-algebra structure on 
$\mbf{U}^+(\mrm C^{\pm})$.

Denote by $\mathcal{V}$  the quotient of the algebra 
$\mathcal{V}'$ by the two-sided ideal generated by the following relations.
\begin{enumerate}
 \item[($\mathcal{V}3$)]
       $E^+_iE^-_i-E^-_iE^+_i=\frac{K_i-K_{-i}}{v_i-v_i^{-1}}\quad$ for any $i\in I$.
\end{enumerate}
In other words, $\mathcal{V}$ is a $\mathbb{Q}(v)$-algebra defined by the following 
generators and relations.
The generators are 
\[
E^+_i, E^-_i, K_i\;\text{ and}\; K_{-i}\quad\text{ for all}\; i\in I.
\]
and the relations are 
\[
(\mathcal{V}1), (\mathcal{V}2)\; \text{ and}\; (\mathcal{V}3)\;
\text{ plus (++)-part, (- -)-part and (+,-)-part.}
\]
Notice that when $i=j$, $(\mathcal{V}2)$ implies that 
\begin{itemize}
 \item[($\star$)] $K_iE^+_i=v_i^2E^+_iK_i$ and $K_iE^-_i=v^{-2}_iE^-_iK_i$, 
                 for any $i\in I$.
\end{itemize}
By the same argument in Lemma ~\ref{sl2}, we see that $(\star)$ and $(\mathcal{V}3)$ imply (a) and 
(b) in (+,-)-part. So $\mathcal{V}$ 
is a $\mathbb{Q}(v)$-algebra defined by the following
generators and relations.
The generators are 
\[
\text{ 
$E^+_i$, $E^-_i$, $K_i$ and $K_{-i}$, for all $i\in I$.
}
\]
and the relations are 
\[
\text{
$(\mathcal{V}1)$, $(\mathcal{V}2)$ and $(\mathcal{V}3)$ plus
                 (++)-part, (- -)-part and (c) in (+\, -)-part.
}
\]
Therefore, we have 
\begin{thm} There is a unique $\mathbb{Q}(v)$-algebra isomorphism 
\[\varphi:\mbf U(\mrm C) \to \mathcal{V},\]
such that $\varphi(E_i)=E^+_i$, $\varphi(F_i)=E^-_i$,  
and $\varphi(K_i)=K_i$ 
for all $i\in I$.
\end{thm}
\begin{proof}
The algebra $\mbf U(\mrm C)$ is a $\mathbb{Q}(v)$-algebra defined by the 
following generators and relations.
\begin{itemize}
 \item Generators: $E_i$, $F_i$, $K_i$ and $K_{-i}$, for all $i\in I$.
 \item Relations: $(\mathcal{V}1)$, $(\mathcal{V}2)$ and $(\mathcal{V}3)$ plus 
       $(++)$-part, $(--)$-part and (c) in $(+\, -)$-part. 
       (Where $E_i^+$ and $E_i^-$ are replaced by $E_i$ and $F_i$, respectively.)
\end{itemize}
So the assignment 
\[E_i\mapsto E^+_i, \quad F_i\mapsto E^-_i \quad \text{and} \quad K_i\mapsto K_i \]
defines  a surjective $\mathbb Q(v)$-algebra homomorphism. Its inverse is given by the assignment
\[
E^+_i\mapsto E_i, \quad E^-_i\mapsto F_i \quad \text{and} \quad K_i\mapsto K_i.
\]
\end{proof}

\section{The application to 
module categories of finite dimensional hereditary algebras}

\subsection{Valued quivers} 
\label{valued}
We recall some definitions and results from [DR].
A valued graph is a pair $\Gamma=(I,\mbf d)$, 
where $I$ is a finite set (of vertices) and 
$\mbf d=(d_{ij})_{i,j\in I}$ is  a set of nonnegative integers satisfying
$d_{ii}=0$ for all $i\in I$ and 
there exist positive integers $d=(d_i)_{i\in I}$
such that $d_{ij}d_j=d_{ji}d_i$ for all $i,j\in I$.

For a given valued graph $\Gamma=(I,\mbf d))$, we can associate to $\Gamma$
a symmetrizable Cartan matrix $\mrm{C}=\mrm{C}(\Gamma)=(a_{ij})_{i\in I}$ by
\[
a_{ii}=2\quad \text{and}\quad  a_{ij}=-d_{ij}\quad  \text{for}\; i\neq j\in I.
\]
The map $\Gamma \mapsto \mrm C(\Gamma)$ then defines 
a one-to-one correspondence between valued graphs and symmetrizable Cartan matrices.

Given any valued graph $(I,\mbf d)$, define $\Omega$ to be the set of all pairs $(i,j)$ in $I$
such that $d_{ij}\neq 0$. By defining two maps $s$ and $t$ from $\Omega$ to $I$ such that 
\[
\{s((i,j)),t((i,j))\}=\{i,j\}\quad \text{for all}\; (i,j)\in \Omega,\] 
the quadruple $Q=(I,\Omega,s,t)$ becomes a quiver. 
We call the pair $(Q, \mbf d)$ a $valued$ $quiver$ associated to $(I,\mbf d)$. 
Note that for a given valued graph, we can associate a lot of valued quivers to it. 

{\em Throughout this paper we assume that the valued quivers are connected and contains 
no oriented cycles.}

Let $k$ be a finite field of $q$ elements and 
$\overrightarrow{\Gamma}=(Q,\mbf d)$ a valued quiver. 
For each $i\in I$, we associate a finite extension $k_i$ of $k$ in an algebraic closure of $k$.
A reduced $k$-species of $\overrightarrow{\Gamma}$ is 
a pair 
\[\mathscr{S}=\left\{(k_i)_{i\in I},\;(_{s(h)}\! M_{t(h)})_{h\in \Omega}\right\}\]
where 
for all $h \in \Omega$, (we write $s(h)=i$, $t(h)=j$), $_iM_j$ is an $k_i$-$k_j$-bimodule such that
\[
\dim \, ( \, _iM_j)_{k_i}=d_{ij}\quad \text{and}\quad  \dim_k \;k_i=d_i.
\]

A representation of $\mathscr{S}$ is a pair 
$\{(V_i)_{i\in I}, (\, _j\varphi_i)_{h\in \Omega}\}$,
($s(h)=i$, $t(h)=j$), where $V_i$ are $k_i$-vector spaces and 
$_j\varphi_i: V_i \otimes \, _iM_j\to V_j$ are $k_j$-linear maps. 
A representation is called finite dimensional if $V_i$ are finite dimensional.
Denote by Rep-$\mathscr S$ the category of 
the finite dimensional representations of $\mathscr S$.

Given a species $\mathscr S$, we define its tensor algebra $\Lambda$ by
\[\Lambda=\otimes_{t\geq 0}\Lambda^{(t)}\]
where (we write $s(h)=i, t(h)=j$)
\[
\Lambda^{(0)}=\prod_{i\in I}k_i,\quad 
\Lambda^{(1)}=\prod_{h\in \Omega}\;_iM_j, \quad
\text{ and} 
\;
\Lambda^{(n)}=\Lambda^{(n-1)}\otimes_{\Lambda^{(0)}}\Lambda^{(1)}\quad \text{for}\; n\geq 2,
\]
with the componentwise addition and the multiplication induced by taking tensor products.
$\Lambda$ is a finite dimensional hereditary $k$-algebra. Furthermore,
any finite dimensional hereditary algebra over $k$ can be obtained in this way.
Note that the category $\Lambda$-mod of finite dimensional left $\Lambda$-modules
is equivalent to the category Rep-$\mathscr S$ of finite dimensional representations of 
$\mathscr S$.
In this paper, we simply identify representations of $\mathscr S$ with
$\Lambda$-modules.

\subsection{Product valued quivers}
Given two valued quivers $\overset{\rightarrow}{\Gamma}=(Q,\mbf{d})$
and 
$\overset{\rightarrow}{\Gamma'}=(Q',\mbf d')$, 
define a valued quiver 
$\overset{\rightarrow}{\Gamma''}=(Q'',\mbf d'')$ 
by
\begin{enumerate}
 \item $I''=I\times I'$.
 \item $d''_{(i,i'),(j,j')}=
       \left\{
        \begin{array}{lll}
         d_{ij} & \mbox{if $i'=j'$},\\
         d'_{i'j'} & \mbox{if $i=j$},\\
         0 & \mbox{otherwise}.
        \end{array}
       \right.$
 \item $\Omega''=\{((i,i'),(j,i'))\,|\, (i,j)\in \Omega, i'\in I'\}
                 \cup
                 \{((i,i'), (i,j'))\,|\, i\in I, (i',j')\in \Omega'\}$
               
                $\;\,=\Omega\times I' \cup I\times \Omega'$ (disjoint union).
 \item $s''((a,b))=
           \left\{
              \begin{array}{ll}
               (s(a),b) & \mbox{if $(a,b)\in \Omega \times I'$},\\
               (a,s'(b)) & \mbox{if $(a,b)\in I\times \Omega'$}.
              \end{array}
           \right.$
 \item $t''((a,b))=
             \left\{
              \begin{array}{ll}
                (t(a), b) & \mbox{if $(a,b)\in \Omega \times I'$},\\
                (a,t'(b)) & \mbox{if $(a,b)\in I\times \Omega'$}.
              \end{array}
              \right.$
\end{enumerate}
(One can check that $\overset{\rightarrow}{\Gamma''}$ is a valued quiver. 
To find $d''=(d''_{(i,i')})_{(i,i')\in I''}$ for $\mbf d''$, 
let $d''_{(i,i')}=d_id'_{i'}$.)
We call $\overset{\rightarrow}{\Gamma''}$ 
the $product$ $valued$ $quiver$ of $\overset{\rightarrow}{\Gamma}$ and
$\overset{\rightarrow}{\Gamma'}$ and  write $\overset{\rightarrow}{\Gamma} \times \overset{\rightarrow}{\Gamma'}$ for it. 

The product valued quiver we study in this paper will be 
$\overset{\rightarrow}{\Gamma^{\pm}}=
(+\overset{(2,2)}{\to}-)\times \overset{\rightarrow}{\Gamma}$, 
for any valued quiver $\overset{\rightarrow}{\Gamma}$.
It can be visualized by taking two copies of
$\overset{\rightarrow}{\Gamma}$,
one positive and one negative,
and connecting the corresponding
vertex in the two copies of 
$\overset{\rightarrow}{\Gamma}$ with an arrow of value $(2,2)$ 
from positive part to the negative part.

Given a species $\mathscr S=(k_i,\,_iM_j)$ of $\overrightarrow{\Gamma}$, we define 
a species $\mathscr S^{\pm}$ of $\overrightarrow{\Gamma}^{\pm}$ by
\begin{enumerate}
\item[(4)] $k_{(+,i)}=k_i=k_{(-,i)}$, for all $i\in I$;
\item[(5)] $_{(+,i)}M_{(+,j)}=\,_iM_j=\,_{(-,i)}M_{(-,j)}$, 
      for all $h\in \Omega$ such that $s(h)=i$ and $t(h)=j$;
\item[(6)] $_{(+,i)}M_{(-,i)}=k_i^2$, where $k_i^2$ is a finite extension field of $k$ such that 
      $\mrm{dim}(k_i^2)_{k_i}=2$. 
\end{enumerate}
One can check that this is a species of $\overrightarrow{\Gamma}^{\pm}$. 
Denote by $\Lambda^{\pm}$ the tensor algebra of $\mathscr S^{\pm}$.
Since $k_{(+,i)}=k_i$ and $_{(+,i)}M_{(+,j)}=\,_iM_j$ for all $i,j\in I$,
there is a natural $k$-algebra embedding
\begin{enumerate}
\item[(7)] $+:\Lambda \to \Lambda^{\pm}$,
\end{enumerate}
such that elements in $\Lambda^0$ (resp. $\Lambda^1$) sends to elements in 
$(\Lambda^{\pm})^0$ (resp. $(\Lambda^{\pm})^1$) of whose
$(-,i)$-components (resp. $(-,h)$-components and $(+\to-,i)$-components) are zero 
for all $i\in I$ (resp. $h\in \Omega$ and $i\in I$).
Similarly, we have a $k$-algebra embedding
\begin{enumerate}
\item[(8)] $-:\Lambda \to \Lambda^{\pm}$,
\end{enumerate}
such that elements in $\Lambda^0$ (resp. $\Lambda^1$) sends to elements in 
$(\Lambda^{\pm})^0$ (resp. $(\Lambda^{\pm})^1$) of whose
$(+,i)$-components (resp. $(+,h)$-components and $(+\to-,i)$-components) are zero 
for all $i\in I$ (resp. $h\in \Omega$ and $i\in I$). 

\subsection{Twisted Ringel-Hall algebras associated to $\Lambda$.}
Let $k$ be a finite field of $q$ elements. Let $v=\sqrt{q}$. 
Let $\Lambda$ be the tensor algebra of a species $\mathscr S$ of 
a given valued quiver $\overrightarrow{\Gamma}=((I,\Omega,s,t),\mbf d)$.
Note that $\Lambda$ is 
a finite dimensional hereditary algebra over $k$.

Denote by $\mathscr P$ the set of isomorphism classes of 
finite dimensional left $\Lambda$-modules.
For any $\alpha \in \mathscr P$, we write $|\alpha|$ for its dimension vector.
Denote by $\{S_i| i\in I\}$ the set of all pairwise non-isomorphic simple $\Lambda$-modules.

The Euler bilinear from 
$<,>\;: \mathbb Z[I]\times \mathbb Z[I] \to \mathbb Z$ is defined by 
\[<i,j>= \mrm{dim}_k\; \mrm{Hom}_{\Lambda}(S_i, S_j) -
\mrm{dim}_k \; \mrm{Ext}^1_{\Lambda}(S_i, S_j)
\quad \text{for any $i, j\in I$}.
\]
The symmetric Euler form $(,)$ is defined by  
\[(a, b)=<a,b>+<a, b> \quad 
\text{
for any $a, b \in \mathbb Z[I]$.
}
\]
Note that $(,)$ is a symmetric bilinear form. 
We then have (see Section ~\ref{valued})
\[
a_{ij}=2\frac{(i,j)}{(i,i)} \quad \text{for}\; i,j\in I.
\] 

Let  $L, M$ and $N$ be the representativers of 
$\alpha, \beta $ and $\gamma \in \mathscr P$.
Denote by $g^{\gamma}_{\alpha,\beta}$ the number of submodules $L'$ of $L$ such that
$L/L'\cong M$ and  $L'\cong N$. Note that $g^{\gamma}_{\alpha,\beta}$ is independent of the choices of the 
representatives of $\alpha,\beta$ and $\gamma$.
 
Let $R$ be an integral domain containing the rational field $\mathbb Q$ and 
the element $v$ and $v^{-1}$.
The $twisted$ $Ringel$-$Hall$ $algebra$ $\mbf{H=H}(\Lambda)$ is 
by definition the free $R$-module
with basis $\{u_{\alpha}|\alpha\in \mathscr P\}$. The multiplication on $\mbf H$ is defined by
\[
u_{\alpha}u_{\beta}=v^{<|\alpha|,|\beta|>}\sum_{\gamma\in\mathscr P}
g^{\gamma}_{\alpha,\beta}\;u_{\gamma}
\quad
\text{for all}\; \alpha,\beta\in \mathscr P.
\]
$\mbf H$ is then an associative ($\mathbb{N}[I]$-graded) $R$-algebra with identity $u_0$ (see ~\cite{Ringel2}).
The grading is defined by 
\[\mbf H=\oplus_{\nu\in\mathbb{N}[I]}\mbf H_{\nu},\] 
where $\mbf H_{\nu}$ is the free $R$-module generated by $u_{\alpha}$ with $\alpha$ of 
dimension vector $\nu$.

For any $\alpha\in\mathscr P$, let $a_{\alpha}$ be the cardinality of
$\mrm{Ext}^1(M,M)$ where $M$ is  a representative of the isomorphism class $\alpha$.
By ~\cite{Green},
there is a symmetric bilinear non-degenerate form 
\[
[\,,\,]: \mbf H\times \mbf H\to R
\] 
defined by
\[
[u_{\alpha},u_{\beta}]=\frac{\delta_{\alpha,\beta}}{a_{\alpha}}\quad \text{for all}\; \alpha,\beta \in \mathscr P.
\]
Given  any $\nu\in \mathbb N[I]$, Let
\[
\mbf L_{\nu}=\sum_{\beta +\gamma=\nu; \beta, \gamma \neq \nu}\mbf H_{\beta}\mbf H_{\gamma}.
\] 
It is an $R$-submodule of $\mbf H_{\nu}$. 
Let 
\[
\mbf L_{\nu}^{\perp}=\{x\in \mbf H_{\nu}|[x, \mbf L_{\nu}]=0\}.
\] 
It is the complement of $\mbf L_{\nu}$ in $\mbf H_{\nu}$.
Since the bilinear form $[,]$ is non-degenerate, 
we can choose an orthonormal basis for each $\mbf L_{\nu}^{\perp}$.
Let $(u_{\mbf i})_{\mbf{i \in I}}$ be the union of these bases. 
Define $|\mbf{i}|=\nu$ if $u_{\mbf i}\in \mbf H_{\nu}$.

Define a bilinear form $(,)_{\mbf I}:\mathbb Z[\mbf I]\times \mathbb Z[\mbf I]\to \mathbb Z$
by
\[ (\mbf{i,j})_{\mbf{I}}=(|\mbf i|, \mbf j|),\;
\mbox{for any} \;\mbf{i,j\in I},\]
where $(,)$ is the Euler form associated to $\Lambda$.
From ~\cite[3.2]{SV1}, the bilinear form $(,)_{\mbf I}$ 
on $\mathbb{Z}[\mbf I]$ satisfies the 
following conditions. 
\begin{enumerate}
\item $(\mbf{i,j})\leq 0$ for any $\mbf{i \neq j \in I}$.
\item $2\frac{(\mbf{i,j})}{(\mbf{i,i})}$ is
      an integer,  if $(\mbf{i,i})> 0$, for all $\mbf{j\in I}$.
\item $(\mbf{i,i})>0$ if and only if $u_\mbf{i}\in \cup_{i\in I} \mbf H_i$.
\item $\mrm{dim} \;  \mbf H_i=1$ for all $i\in I$.
\end{enumerate}
When $(\mbf{i,i})>0$, we can choose $u_{\mbf{i}}$ to be the isomorphism class of 
the simple $\Lambda$-module in $\mbf H_{i_0}$, for some $i_0\in I$. 
By abuse of the notation, we simply write $u_{i_0}$ for it.
So $\mbf I=I \cup \mbf J$ for some suitable index set $\mbf J$. 

If $\mbf i\in I$, let  
\[
a_{\mbf{ij}}=2(\mbf{i,j})_{\mbf I}/(\mbf{i,i})_{\mbf I}.
\] 
If $\mbf i\in \mbf J$, let 
\[a_{\mbf{ij}}=(\mbf{i,j})_{\mbf I}.\] 
The matrix $\mrm{C}_{\mbf I}=(a_{\mbf{ij}})_{\mbf{i,j\in I}}$ is 
then a symmetrizable Borcherds-Cartan matrix.

\subsection{Double Ringel-Hall algebras associated to $\Lambda$.}
Let 
$K_{\mbf i}=K_{|i|}$, $K_{-\mbf i}=K_{-|i|}$ and $v_{\mbf i}=v$
for all $\mbf i\in \mbf J$,
$v_{i}=v^{(i,i)/2}$ for any $i\in I$.

The double Ringel-Hall algebra $\mathscr D(\Lambda)$ (see ~\cite{SV2} and ~\cite{DX1}) 
associated to $\Lambda$ is 
an $R$-algebra characterized by  generators and relations. The generators are
\begin{enumerate}  
       \item[(G1)] $X_{\mbf i}$ and  $Y_{\mbf{i}}\quad$ 
             for all $\mbf{i\in I}$,
       \item[(G2)] $K_i$ and $K_{-i}\quad$ for all $i \in I$, 
\end{enumerate}
and  the  relations are 
      \begin{enumerate} 
      \item[(R1)] $K_0=1$, $K_iK_{-i}=1=K_{- i}K_i$, 
           and  $K_iK_j=K_{i+j}\quad$ for all $i,j\in I$.
      \item[(R2)] $K_i X_{\mbf j}=
            v^{(i,\mbf j)_{\mbf I}} X_{\mbf j} K_i$ 
           and  $K_i Y_{\mbf j}=
            v^{-(i,\mbf j)_{\mbf I}} Y_{\mbf j}K_i\;\;$ 
            for all $i\in I$ and  $\mbf{j\in I}$.
      \item[(R3)] $\sum_{p=0}^{1-a_{i\mbf j}}(-1)^p
            \bin{1-a_{i\mbf j}}{ p}_i
            X_i^p X_{\mbf j} X_i^{1-a_{i\mbf j}-p}=0\quad $ and 
      \item[] $\sum_{p=0}^{1-a_{i\mbf j}} (-1)^p
            \bin{1-a_{i\mbf j}}{ p}_i
            Y_i^p Y_{\mbf j} Y_i^{1-a_{ij}-p}=0\;\;$ 
            for $i\in I$ and $\mbf{j\in I}$ ($\mbf j \neq i$).
      \item[(R4)] $X_{\mbf i}X_{\mbf j}=X_{\mbf j}X_{\mbf i}\;$ and 
            $Y_{\mbf i}Y_{\mbf j}=Y_{\mbf j}Y_{\mbf i}\;$ 
            for any $\mbf{i,j\in I}$ such that $(\mbf{i,j})_{\mbf I}=0$.
      \item[(R5)] $X_{\mbf i}Y_{\mbf j}-Y_{\mbf j} X_{\mbf i}=
             \delta_{\mbf{i,j}}
            \frac{ K_{\mbf i}-K_{-\mbf i}}{v_{\mbf i}-v_{\mbf i}^{-1}}\;\;$
            for all $\mbf{i,j\in I}$. 
     \end{enumerate}
$\mathscr D(\Lambda)$ has a triangular decomposition 
$\mathscr D^+(\Lambda)\otimes \mathscr D^0(\Lambda) 
\otimes \mathscr D^-(\Lambda)$,
where $\mathscr D^+(\Lambda)$ (resp. $\mathscr D^-(\Lambda)$)
is generated by $X_{\mbf i}$ (resp. $Y_{\mbf i}$), 
for all $\mbf{i\in I}$
and $\mathscr D^0(\Lambda)$ is generated by $K_i$ and $K_{-i}$,
for all $i\in I$.
Moreover, there exists a unique $R$-algebra isomorphism 
\[
\tag{$\star$}
 \mathscr D^+(\Lambda) \to \mbf H(\Lambda),
\]
such that $X_{\mbf i} \mapsto u_{\mbf i}$, for all $\mbf{i\in I}$.

\subsection{Extended twisted Ringel-Hall algebras associated to $\Lambda^{\pm}$.}
\label{extended}
We preserve the setting in Section 3.3.
Recall that $\overrightarrow{\Gamma}=(Q,\mbf d)$, 
(with $Q=(I,\Omega,s,t)$), is a valued quiver. 
$\Lambda$ is the tensor algebra
of a $k$-species $\mathscr S$ of $\overrightarrow{\Gamma}$.
We have from Section 3.2 
the tensor algebra $\Lambda^{\pm}$ of the $k$-species $\mathscr S^{\pm}$ 
of $\overrightarrow{\Gamma}^{\pm}$.
Denote by $\mbf{H}(\Lambda^{\pm})$ the twisted Ringel-Hall algebra of $\Lambda^{\pm}$.

By the construction in Section 3.3, 
we have a generating set $(\theta_a)_{a\in A}$ ($A$ an index set) for 
$\mbf H(\Lambda^{\pm})$. We can define a bilinear form 
$(,)_A:\mathbb Z[A]\times\mathbb Z[A]\to \mathbb Z$ as the bilinear form $(,)_{\mbf I}$ in
Section 3.3.
By Section 3.4 ($\star$), there is a unique $R$-algebra
isomorphism 
\[
\tag{1} 
 \mathscr D^+(\Lambda^{\pm}) \to \mbf H(\Lambda^{\pm}),
\]
such that $X_a \mapsto \theta_a$, for all $a\in A$.

Denote by $\mathscr P(\Lambda^{\pm})$ 
the set of isomorphism classes of left $\Lambda^{\pm}$-modules.
Denote by $\mathscr P^+$ (resp. $\mathscr P^-$)  the subset of 
$\mathscr P(\Lambda^{\pm})$
consisting of all isomorphism classes of left $\Lambda^{\pm}$-modules 
of dimension vectors in $\mathbb{Z}[I^+]$
(resp. $\mathbb{Z}[I^-]$).
Denote by $\mbf H^{\pm}$ (resp. $\mbf H^+$, $\mbf H^-$) 
the subalgebra of $\mbf{H}(\Lambda^{\pm})$ generated by
$\mathscr P^+$ and $\mathscr P^-$ (resp. $\mathscr P^+$, $\mathscr
P^-$). 
The two $k$-algebra embeddings 
$+:\Lambda \to \Lambda^{\pm}$ and  $-:\Lambda \to  \Lambda^{\pm}$
in Section 3.2
induce two bijective maps 
\[
+:\mathscr P \to \mathscr P^+\quad
\text{  and }
\quad
-:\mathscr P \to \mathscr P^-,
\]
respectively. We have 

\begin{lem} 
\label{+-}
The maps $+:\mathscr P \to \mathscr P^+$ and 
$-:\mathscr P \to \mathscr P^-$
induce $R$-algebra isomorphisms $+:\mbf H(\Lambda) \to \mbf H^+$ and  
$-:\mbf H(\Lambda) \to \mbf H^-$, 
respectively.
\end{lem}
\begin{proof}
We show that $+$-part holds. The proof of  the $-$-part is similar to the $+$-part. 
Note that $\mbf H(\Lambda)$ and $\mbf H^+$ are isomorphic as 
$R$-modules. We only need to prove that 
$+:\mbf H(\Lambda)\to \mbf H^+$ is an $R$-algebra homomorphism.
It reduces to check that $<\alpha,\beta>=<\alpha^+,\beta^+>$ and 
$g^{\gamma}_{\alpha,\beta}=g^{\gamma^+}_{\alpha^+,\beta^+}$ 
for any $\alpha, \beta, \gamma \in \mathscr P$.
But this can be checked directly from the definitions.
\end{proof}

Given any $\alpha\in \mathscr P$, 
we write $\alpha^+$ (resp. $\alpha^-$)
for the corresponding element $+(\alpha)$ (resp. $- (\alpha))$ in
$\mathscr P^+$ (resp. $\mathscr P^-$).
We also write $x^+$ (resp. $x^-$) for $+(x)$ (resp. $-(x)$) in $\mbf H^+$ (resp. $\mbf H^-$),
for all $x$ in $\mbf H(\Lambda)$. Here $+(\;)$ and $-(\;)$ are the maps in
Lemma ~\ref{+-}.

By Section 3.3, we have a generating set $(u_{\mbf i})_{\mbf{i\in I}}$ for $\mbf H(\Lambda)$. 
By Lemma ~\ref{+-}, we have that its images 
($u_{\mbf{i}}^+)_{\mbf{i\in I}}$ and ($u_{\mbf i}^-)_{\mbf{i \in I}}$
under
$+(\;)$ and $-(\;)$
are generating sets for $\mbf H^+$ and $\mbf H^-$, respectively.

Note that $\mbf H^{\pm}$ is generated by $\mbf H^+$ and $\mbf H^-$, 
so that the set
$\{u^+_{\mbf i}, u^-_{\mbf i}\}_{\mbf{i\in I}}$ is a generating set for
$\mbf H^{\pm}$. By the construction in Section 3.3, 
we can choose the generating set $(\theta_a)_{a\in A}$
for $\mbf H(\Lambda^{\pm})$ 
such that it contains  $\{u^+_{\mbf i}, u^-_{\mbf i}\}_{\mbf{i\in I}}$.
By 3.5 (1), $\mbf H^{\pm}$ is an $R$-algebra 
defined by the following generators and relations.
\begin{enumerate}
\item Generators: $u^+_{\mbf i}$ and  $u^-_{\mbf i},\;\;\;$ for all $\mbf{i\in I}$.
\item Relations:
      \begin{itemize}
       \item [($1^+$)] 
               $\sum_{p=0}^{1-a_{i\mbf j}} (-1)^p
               \bin{1-a_{i\mbf j}}{p}_i   
               (u^+_i)^p\, u^+_{\mbf j} \,(u^+_i)^{1-a_{i\mbf j}-p}=0\;\;$ 
                  
              for any $i\in I$ and $\mbf{j\in I}$ such that $\mbf j \neq i$.
       \item [ ($2^+$)] 
              $u^+_{\mbf i}\,u^+_{\mbf j}=u^+_{\mbf j}\,u^+_{\mbf i}\;$ 
              for any $\mbf{i,j\in I}$ such that
              $(\mbf{i,j})=0$.
       \item [($1^-$)] 
               $\sum_{p=0}^{1-a_{i\mbf j}} (-1)^p
               \bin{1-a_{i\mbf j}}{p}_i
               (u^-_i)^p\, u^-_{\mbf j}\, (u^-_i)^{1-a_{i \mbf j}-p}=0\;\;$ 
                  
               for any $i\in I$ and $\mbf{j\in I}$ such that $\mbf j \neq i$.
       \item [($2^-$)] 
              $u^-_{\mbf i}\,u^-_{\mbf j}=u^-_{\mbf j}\,u^-_{\mbf i}\;$ 
              for any $\mbf{i,j\in I}$ such that
              $(\mbf{i,j})=0$.
       \item [($1^{\pm}$)] 
               $\sum_{p=0}^3 (-1)^p
                \bin{3}{p}_i
                (u^+_i)^p \,u^-_i \,(u^+_i)^{3-p}=0 \;\;$ and 
          
         \item[]      $\sum_{p=0}^3 (-1)^p
                \bin{3}{p}_i
                (u^-_i)^p\, u^+_i\, (u^-_i)^{3-p}=0 \;\;$ 
                for any $i\in I$.
       \item [($2^{\pm}$)] 
               $\sum_{p=0}^{1+b_{i\mbf j}} (-1)^p
                \bin{1+b_{i\mbf j}}{p}_i
                (u^+_i)^p \,u^-_{\mbf j} \,(u^+_i)^{1+b_{i \mbf j}-p}=0\quad $ and 

               $\sum_{p=0}^{1+b_{i\mbf j}} (-1)^p
                \bin{1+b_{i\mbf j}}{p}_i
                 (u^-_i)^p \,u^+_{\mbf j} \,(u^-_i)^{1+b_{i \mbf j}-p}=0\;\;$ 
                  
                  for any $i\in I$ and $\mbf{j\in J}$ ($\mbf j \neq i$) 
                  such that $(i^+,\mbf j^-)_{\Gamma^{\pm}}\neq 0$.
       \item [($3^{\pm}$).] 
              $u^+_{\mbf i}\,u^-_{\mbf j}=u^-_{\mbf j}\,u^+_{\mbf i}\quad $ and 
              $u^-_{\mbf i}\,u^+_{\mbf j}=u^+_{\mbf j}\,u^-_{\mbf i}\;$ 
              for any $\mbf{i,j\in I}$ such that
              $(\mbf{i^+,j^-})_{\Gamma^{\pm}}=0$. 
      \end{itemize}
\end{enumerate}
where 
$b_{i\mbf j}:=4 d_i |\mbf j|_{i}$ for any $i\in I$, and $\mbf{j\in J}$.
($ |\mbf j|=\sum |\mbf j|_i\;i \in \mathbb N[I]$.)
Note that $-b_{i\mbf j}$ is the number $(i^+,\mbf j^-)_{\Gamma^{\pm}}$, 
where $(,)_{\Gamma^{\pm}}$ is 
the symmetric Euler form associated to the valued graph $\Gamma^{\pm}$. 

Consider the $R$-module $\mathscr H'=R[K_i,K_{-i}]_{i\in I}\otimes \mbf H^{\pm}$. 
There is a unique $R$-algebra structure on $\mathscr H'$ such that
\begin{itemize}
\item [$(1^0)$]
      $K_0=1$, $K_iK_{-i}=1=K_{-i}K_i$, and $K_iK_j=K_{i+j}\;\;$ for all $i,j\in I$,
\item [$(2^0)$]
      $K_ix^+=v^{(i,\nu)}x^+K_i$ and 
      $K_ix^-=v^{-(i,\nu)}x^-K_i\quad $
      for all $i\in I$ and  $x\in \mbf H(\Lambda)_{\nu}$,
\end{itemize}
(where $(,)$ is the Euler form associated to $\Lambda$) 
and preserves the $R$-algebra structure on $\mbf H^{\pm}$.

Denote by $\mathscr H$ the quotient of the algebra $\mathscr H'$ by the two-sided ideal
generated by the following relations.
\begin{itemize}
\item [($4^{\pm}$)]
      $u^+_i\,u^-_i-u^-_i \, u^+_i=\frac{K_i-K_{-i}}{v_i-v_i^{-1}}\;\;$ for all $i\in I$.
\item [($5^{\pm}$)]
      $u^+_{\mbf j}\,u^-_{\mbf j}-u^-_{\mbf j} \, u^+_{\mbf j}
       =\frac{K_{\mbf j}-K_{-\mbf j}}{v-v^{-1}}\;\;$ for all $\mbf{ j\in J}$.
\item [($6^{\pm}$)]
      $u^+_i\,u^-_{\mbf j}=u^-_{\mbf j} \, u^+_i\;\;$ and  
      $u^-_i\,u^+_{\mbf j}=u^+_{\mbf j} \, u^-_i$ 
      for all $i\in I$ and $\mbf{j\in J}$ 
      such that $(i^+,\mbf j^-)_{\Gamma^{\pm}}\neq 0$.
\end{itemize}
We have 
\begin{thm}
\label{main}
There exists a unique $R$-algebra isomorphism 
\[\phi:\mathscr D(\Lambda) \to \mathscr H,\]
such that $\phi(X_{\mbf i})=u^+_{\mbf i}$, 
$\phi(Y_{\mbf i})=u^-_{\mbf i}$ and    
$\phi(K_i)=K_i$, for all $\mbf{i\in I}$ and $i\in I$.
\end{thm}
The proof will be given in the next subsection.

Let $\mathscr C$ be the subalgebra of $\mathscr H$ generated by $u^+_i$, $u^-_i$, 
$K_i$ and $K_{-i}$, for all
$i\in I$. Let $\mbf U(\mrm C)$ be 
the quantized enveloping algebra associated to $\mrm C$. 
(See Section 2.1) 
Here $\mrm{C}=(a_{ij})_{i\in I}$ is  
the symmetrizable generalized Cartan matrix associated to 
$\Lambda$.
Let $\mbf{U}_q(\mrm C)$ be its specialization at $v=\sqrt{q}$.
A consequence of Theorem ~\ref{main} is as follows.

\begin{cor} There is a unique $R$-algebra isomorphism 
$\varphi: \mbf{U}_q(\mrm C) \to \mathscr{C}$
such that $\varphi(E_i)=u^+_i$, $\varphi(F_i)=u^-_i$ and $\varphi(K_i)=K_i$ for all $i\in I$.
\end{cor}
\begin{proof} 
In fact, the Cartan matrix associated to the underlying graph of 
$\overrightarrow{\Gamma}^{\pm}$ is nothing
but the Cartan matrix $\mrm{C}^{\pm}$ constructed in Section 2.1. 
So the Corollary is just a $v=\sqrt{q}$ version of Theorem 2.5.
\end{proof} 

{\bf Remark.} We only discuss the realization of the specialization $\mathbf U_q(\mrm C)$ of the quantum group $\mbf U$ 
via the representations of $\Lambda^{\pm}$. In fact, one can realize $\mbf U$ via the generic composition algebra
of $\Lambda^{\pm}$ by the above process essentially word by word.

\subsection{Proof of Theorem 3.7.}
By definition, we see that $\mathscr H$ is an $R$-algebra defined by the following generators
and relations.
The generators are 
\[
\text{ 
$K_{i}$, $K_{-i}$, $u^+_{\mbf i}$ and $u^-_{\mbf{i}}\;\;$, for all $i\in I$ and $\mbf{i\in I}$.
}
\] 
and the relations are 
\[
\text{
$(1^0)$, $(2^0)$, $(1^+)$, $(2^+)$, $(1^-)$, $(2^-)$, $(1^{\pm})$, $(2^{\pm})$
$(3^{\pm})$, $(4^{\pm})$, $(5^{\pm})$ and $(6^{\pm})$ in Section ~\ref{extended}. 
}
\]
On the other hand, the Ringel-Hall algebra $\mathscr D(\Lambda)$ 
is an $R$-algebra defined by the generators and
relations listed in Section 3.4 (G1-G2) and (R1-R5).
By comparing the generators and the relations, we see that the 
relations  $(1^0)$, $(2^0)$, $(1^+)$, $(2^+)$, $(1^-)$, $(2^-)$,
           $(3^{\pm})$, $(4^{\pm})$, $(5^{\pm})$ and $(6^{\pm})$ 
are exactly all relations $\mathscr D(\Lambda)$ satisfied. 
(With $u^+_{\mbf i}$ and $u^-_{\mbf{i}}$ replaced
by $X_{\mbf i}$ and $Y_{\mbf i}$, respectively, for all $\mbf{i \in I}$.)
By the universality of $\mathscr D(\Lambda)$, 
we have a unique $R$-algebra homomorphism
\[\phi:\mathscr D(\Lambda) \to \mathscr H,\]
such that $\phi(X_{\mbf i})=u^+_{\mbf i}$, 
$\phi(Y_{\mbf i})=u^-_{\mbf i}$,    
$\phi(K_i)=K_i$, for all $\mbf{i\in I}$ and $i\in I$.
$\phi$ is surjective since $X_{\mbf i}, Y_{\mbf i}$ and $K_i$ ($\mbf{i\in I}$ and $i\in I$) are generators of $\mathscr H$. 

To show that $\phi$ is injective, it suffices to prove the following two statements.
\begin{itemize}
\item [(S1)]
      The relations $(2^0)$, $(4^{\pm})$, and $(5^{\pm})$ imply the relations 
      $(1^{\pm})$.
\item [(S2)]
      The relation $(6^{\pm})$ implies the relation $(2^{\pm})$. 
\end{itemize} 
Note that (S1) has been proved in the proof of Lemma 2.3 with $v$ replaced by $v_i$.
We only need to prove (S2). In $(2^{\pm})$, 
$b_{i,\mbf j}:=4 d_i |\mbf j|_{i}$ is an even number, we write
$b_{i,\mbf j}=2n$ for some $n$. For simplicity, we write $N=2n+1$.
We simplify the left hand side of $(2^{\pm})$ as follows.
\begin{itemize}
 \item[LHS=] $\sum_{p=0}^{N} (-1)^p
               \bin{N}{p}_i
                 (u^+_i)^p \,u^-_{\mbf j} \,(u^+_i)^{N-p},$
 \item[=]    $\sum_{p=0}^{n} \left \{(-1)^p
                \bin{N}{p}_i
                 (u^+_i)^p \,u^-_{\mbf j}\, (u^+_i)^{N-p}
                 + (-1)^{N-p}
                 \bin{N}{N-p}_i
                 (u^+_i)^{N-p}\, u^-_{\mbf j} \,(u^+_i)^{p}\right\}$,
\item[=]      $\sum_{p=0}^{n} (-1)^p 
                \bin{N}{p}_i
                 \left\{(u^+_i)^p \,u^-_{\mbf j} \,(u^+_i)^{N-p}-
                 (u^+_i)^{N-p}\, u^-_{\mbf j}\, (u^+_i)^{p}\right\}$,

\item[=]      $\sum_{p=0}^{n} (-1)^{p+1} 
                \bin{N}{p}_i
                (u^+_i)^p\,\left\{(u^+_i)^{N-2p}\,
                u^-_{\mbf j}-u^-_{\mbf j}\, (u^+_i)^{N-2p}
                \right\}\, (u^+_i)^p$.
\end{itemize}
By the fact that given any positive integer $m$, for any $x$ and $y$,  
\[x^my-yx^m=\sum_{a+a'=m-1}x^a(xy-yx)x^{a'}.\]
We have Statement (S2). Theorem 3.7 follows.

\section{Generalization of the constructions}
For a given  Cartan matrix $\mrm C$, one can give a realization of the quantum group associated to 
$\mrm C$ by starting from
the Cartan matrix $\mrm C^{2n}$, ($n$ any positive integer), instead of $\mrm C^{\pm}$ as in
Theorem 2.5. Here
\[
\mrm C^{2n}=\left [
          \begin{array}{cc}
          \mrm C & -2n\,\mrm{Id}\\
          -2n\,\mrm{Id} & \mrm C
          \end{array}
         \right ],
\]
where $\mrm{Id}$ is an identity matrix with suitable rank.
Note that when $n=1$, $\mrm C^{2n}$ is just $\mrm C^{\pm}$. 
One can also generalize Theorem 3.7 to the case when $\overrightarrow{\Gamma}^{\pm}$ is
replaced by $\overrightarrow{\Gamma}^{2n}=
(+\overset{(2n,2n)}{\to}-)\times \overrightarrow{\Gamma}$ ($n$ any positive integer).
The crucial point in the process of the generalization  is the following lemma, 
which is more or less a generalization of Lemma 2.3.
\begin{lem}
The relations
\begin{itemize}
\item [(1)] $KE^+=v^{2}E^+K$, $KE^-=v^{-2}E^-K$,
\item [(2)] $E^+E^--E^-E^+=\frac{K-K^{-1}}{v-v^{-1}}$.
\end{itemize}
implies the relations
\begin{itemize}
\item [(3)]$\sum_{p=0}^{2n+1}(-1)^p
             \bin{2n+1}{p}              
             (E^+)^p E^- (E^+)^{2n+1-p}=0,$
\item [(4)]$\sum_{p=0}^{2n+1}(-1)^p
              \bin{2n+1}{p}   
              (E^-)^p E^+ (E^-)^{2n+1-p}=0.$ 
\end{itemize}
\end{lem} 
\begin{proof}
We only prove that the relations (1) and (2) imply the relation (3).
(4) can be proven similarly.
As in the proof of Theorem3.7, the left hand side of (3) can be written as follows.

\begin{itemize}
\item[LHS=]$\sum_{p=0}^{2n+1}(-1)^p
              \left[
                \begin{array}{c}
                   2n+1\\ p
                \end{array}
              \right ]
                (E^+)^p E^- (E^+)^{2n+1-p},$
\item[=]$\sum_{p=0}^{n} (-1)^{p+1} \bin{2n+1}{p}
                (E^+)^p \,\left((E^+)^{2n+1-2p}E^- - E^-(E^+)^{2n+1-2p}\right) \,
                (E^+)^p.$
\end{itemize}
By the fact that 
\[(E^+)^mE^--E^-(E^{+})^m=\sum_{a=0}^{m-1}(E^+)^a(E^+E^--E^-E^+)(E^+)^{m-1-a}.\]
We have 
\begin{itemize}
\item [LHS=]
$\sum_{p=0}^{n} (-1)^{p+1} \bin{2n+1}{p}
                (E^+)^p \,\left\{
               \sum_{a=0}^{2n-2p}(E^+)^a(E^+E^--E^-E^+)(E^+)^{2n-2p-a}\right\}(E^+)^p,$
\item [=]
$\sum_{p=0}^{n} (-1)^{p+1} \bin{2n+1}{p}
                \sum_{a=0}^{2n-2p}(E^+)^{p+a}(E^+E^--E^-E^+)(E^+)^{2n-p-a}.$
\end{itemize}
By (2), we have
\begin{itemize}
\item[LHS=] $\sum_{p=0}^{n} (-1)^{p+1} \bin{2n+1}{p}
                \sum_{a=0}^{2n-2p}(E^+)^{p+a}\;\frac{K-K^{-1}}{v-v^{-1}}\;(E^+)^{2n-p-a},$

\item[=] $\frac{1}{v-v^{-1}}\sum_{p=0}^{n} (-1)^{p+1} \bin{2n+1}{p}
                \sum_{a=0}^{2n-2p}(E^+)^{p+a}\;K\;(E^+)^{2n-p-a}$
\item[]  $-\frac{1}{v-v^{-1}}\sum_{p=0}^{n} (-1)^{p+1} \bin{2n+1}{p}
                \sum_{a=0}^{2n-2p}(E^+)^{p+a}\;K^{-1}\;(E^+)^{2n-p-a}.$
\end{itemize}
It suffices to prove that the first term is zero. Note that 
\begin{itemize}
\item[(A)] $ \sum_{a=0}^{2n-2p}(E^+)^{p+a}\;K\;(E^+)^{2n-p-a}$
\item [=] $\sum_{a=0}^{n-p}(E^+)^{p+a}\;K\;(E^+)^{2n-p-a}+
          \sum_{a=n-p+1}^{2n-2p}(E^+)^{p+a}\;K\;(E^+)^{2n-p-a},$
\item [$\overset{(\star)}{=}$]
          $\sum_{a=0}^{n-p}(v^2)^{n-(p+a)}(E^+)^n\;K\;(E^+)^n+
          \sum_{a=n-p+1}^{2n-2p}(v^{-2})^{(p+a)-n}(E^+)^n\;K\;(E^+)^n,$
\item [=] $\left\{\sum_{a=0}^{n-p}(v^2)^{n-(p+a)}+
          \sum_{a=n-p+1}^{2n-2p}(v^{-2})^{(p+a)-n}\right\}
          (E^+)^n\;K\;(E^+)^n,$
\item[=]  $\left\{\sum_{a=0}^{n-p}(v^2)^{(n-p)+a}+
          \sum_{a=n-p+1}^{2n-2p}(v^{2})^{(n-p)-a}\right\}
          (E^+)^n\;K\;(E^+)^n,$
\item[=]  $\left\{\sum_{a=0}^{n-p}(v^2)^a+
          \sum_{a=-(n-p)}^{-1}(v^{2})^a\right\}
          (E^+)^n\;K\;(E^+)^n,$
\item[=]  $\sum_{a=-(n-p)}^{n-p}v^{2a}
          (E^+)^n\;K\;(E^+)^n,$
\item[=]  $[2(n-p)+1](E^+)^n\;K\;(E^+)^n.$
\end{itemize}
Here the equation $(\star)$ holds by applying the relation (1). 
By (A), to prove the Lemma, it suffices to show that the following is zero:
\begin{itemize}
\item [(B)] $\sum_{p=0}^{n} (-1)^{p+1} \bin{2n+1}{p}
                 [2(n-p)+1].$
\end{itemize}
Let $B_{n,i}=\sum_{p=n-i+1}^{n} (-1)^{p+1} \bin{2n+1}{p}
                 [2(n-p)+1]$.
We have
\begin{lem}
$B_{n,i}=-(-1)^{n-i}\bin{2n+1}{n-i}
                   \frac{[n+i+1][i]}{[n]},$
for all $1\leq i \leq n$.
\end{lem}
In particular, $B_{n,n}=-[2n+1]$.
By definitions, (B)= $[2n+1]+B_{n,n}$. Therefore, (B)= $0$. 
This finishes the proof of Lemma 4.1.
\end{proof}
Proof of Lemma 4.2. We prove by induction.
When $i=1$, it holds trivially. Now
\begin{itemize}
\item [$B_{n,i+1}$] $=B_{n,i}+(-1)^{n-i} \bin{2n+1}{n-i}[2i+1]$,
\item [$\overset{(\star)}{=}$]
                     $(-1)^{n-i} \bin{2n+1}{n-i}
                      \left([2i+1]- \frac{[n+i+1][i]}{[n]} \right)$,
\item [=]           $(-1)^{n-i} \bin{2n+1}{n-i}
                       \frac{[2i+1][n]-[n+i+1][i]}{[n]}$,
\item [$\overset{(\star\star)}{=}$]             
                     $(-1)^{n-i} \bin{2n+1}{n-i}
                        \frac{[n-i][i+1]}{[n]}$
\item [=] $-(-1)^{n-(i+1)} \bin{2n+1}{n-(i+1)}
                   \frac{[n+(i+1)+1][i+1]}{[n]}.$
\end{itemize}
Here the equality $(\star)$ is by induction. $(\star\star)$ is by the fact that
\[
[2i+1][n]-[n+i+1][i]=[n-i][i+1]\quad \text{ for $1 \leq i \leq n$},
\]
which can be proved by definition. Lemma 4.2 follows.

\end{document}